\documentclass{ifacconf}

\usepackage{graphicx,amsmath,amsfonts,amssymb,mathrsfs,bm}
\usepackage[dvipsnames]{xcolor}
\usepackage{natbib}

\newcommand{\differential}{{\rm{d}}}

\newcommand{\prox}{{\rm{prox}}}

\newtheorem{proposition}{Proposition}

\newtheorem{theorem}{Theorem}

\begin{document}
\begin{frontmatter}

\title{A Distributed Algorithm for Measure-valued Optimization with Additive Objective\thanksref{footnoteinfo}} 

\thanks[footnoteinfo]{This work is partially supported by NSF grants 1923278, 2112755.}

\author[First]{Iman Nodozi} 
\author[Second]{Abhishek Halder} 

\address[First]{Department of Electrical and Computer Engineering, University of California, Santa Cruz, CA 95064, USA (e-mail: inodozi@ucsc.edu).}
\address[Second]{Department of Applied Mathematics, University of California, Santa Cruz, CA 95064, USA (e-mail: ahalder@ucsc.edu)}

\begin{abstract}                
We propose a distributed nonparametric algorithm for solving measure-valued optimization problems with additive objectives. Such problems arise in several contexts in stochastic learning and control including Langevin sampling from an unnormalized prior, mean field neural network learning and Wasserstein gradient flows. The proposed algorithm comprises a two-layer alternating direction method of multipliers (ADMM). The outer-layer ADMM generalizes the Euclidean consensus ADMM to the Wasserstein consensus ADMM, and to its entropy-regularized version Sinkhorn consensus ADMM. The inner-layer ADMM turns out to be a specific instance of the standard Euclidean ADMM. The overall algorithm realizes operator splitting for gradient flows in the manifold of probability measures. 
\end{abstract}


\end{frontmatter}

\section{Introduction}
We consider measure-valued optimization problems of the form 
\begin{align}
\underset{\mu\in\mathcal{P}_{2}(\mathbb{R}^{d})}{\arg\inf} F_{1}(\mu) + F_{2}(\mu) + \hdots + F_{n}(\mu)
\label{AdditiveOptimizationMeasure}    
\end{align}
for some finite integer $n>1$, where $\mathcal{P}_{2}(\mathbb{R}^{d})$ denotes the space of Borel probability measures over $\mathbb{R}^{d}$ with finite second moments. We suppose that the functionals $F_{i}:\mathcal{P}_{2}(\mathbb{R}^{d})\mapsto \mathbb{R}$ are convex for all $i\in[n]$. If the optimization in $\eqref{AdditiveOptimizationMeasure}$ is instead over $\mathcal{P}_{2,{\rm{ac}}}(\mathbb{R}^{d})$, defined as the subset of $\mathcal{P}_{2}(\mathbb{R}^{d})$ such that its elements are absolutely continuous w.r.t. the Lebesgue
measure, then we can rewrite\footnote{with slight abuse of notation in the sense \eqref{AdditiveOptimizationPDF} uses the same symbols $F_{i}$ as in \eqref{AdditiveOptimizationMeasure} for the additive functionals.} \eqref{AdditiveOptimizationMeasure} as 
\begin{align}
\underset{\rho}{\arg\inf} \: F_{1}(\rho) + F_{2}(\rho) + \hdots + F_{n}(\rho)
\label{AdditiveOptimizationPDF}    
\end{align}
where the decision variable $\rho$ is a joint PDF over $\mathbb{R}^{d}$ with finite second moment. 

Problems of the form \eqref{AdditiveOptimizationMeasure} and \eqref{AdditiveOptimizationPDF} arise in several contexts in statistics, machine learning, and control theory. This includes sampling from an unnormalized prior via Langevin Monte Carlo (see e.g., \cite{stramer1999langevin1,stramer1999langevin2,jarner2000geometric,roberts2002langevin,vempala2019rapid}), policy optimization in reinforcement learning (see e.g., \cite{zhang2018policy,chu2019probability,zhang2020variational}), stochastic prediction (see e.g., \cite{jordan1998variational,ambrosio2005gradient,caluya2019proximal,caluya2019gradient}) and estimation (see e.g., \cite{halder2017gradient,halder2018gradient,halder2019proximal}), density control (see e.g., \cite{caluya2021wasserstein,caluya2021reflected}), mean field analysis of neural supervised (see e.g., \cite{chizat2018global,mei2018mean,rotskoff2018neural,sirignano2020mean}) and unsupervised learning (see e.g., \cite{domingo2020mean}).

Let $F := F_{1} + \hdots + F_{n}$. There is a natural connection between problems of the form \eqref{AdditiveOptimizationMeasure} and that of the Wasserstein gradient flow
{\small{\begin{align}
\dfrac{\partial\mu}{\partial t} = - \nabla^{W_{2}}F(\mu) := \nabla\cdot\left(\mu\dfrac{\delta F}{\delta\mu}\right),   
\label{WassgradFlow}    
\end{align}}}
where $\nabla$ denotes the $d$ dimensional Euclidean gradient, and $\frac{\delta}{\delta\mu}$ denotes the functional derivative w.r.t. $\mu$. The operator $\nabla^{W_{2}}$ in \eqref{WassgradFlow} denotes the gradient w.r.t. the 2-Wasserstein metric $W_{2}$ between a pair of probability measures $\mu_{x},\mu_{y}\in\mathcal{P}_{2}\left(\mathbb{R}^{d}\right)$, defined as
{\small{\begin{align}
W_{2}\left(\mu_x,\mu_y\right) := \left(\underset{\pi\in\Pi\left(\mu_x,\mu_y\right)}{\inf}\displaystyle\int_{\mathbb{R}^{2d}}\!\!c\left(\bm{x},\bm{y}\right)\:\differential\pi(\bm{x},\bm{y})\right)^{\frac{1}{2}},
\label{DefWassContinuous}    
\end{align}}}
where $\Pi\left(\mu_x,\mu_y\right)$ is the set of joint probability measures or couplings over the product space $\mathbb{R}^{2d}$, having $\bm{x}$ marginal $\mu_x$, and $\bm{y}$ marginal $\mu_{y}$. We use the ground cost $c\left(\bm{x},\bm{y}\right):=\|\bm{x}-\bm{y}\|_{2}^{2}$, the squared Euclidean distance in $\mathbb{R}^{d}$. It is well-known \cite[Ch. 7]{villani2003topics} that $W_{2}$ defines a metric on $\mathcal{P}_{2}\left(\mathbb{R}^{d}\right)$. For notational ease, we henceforth drop the subscript from $W_{2}$, and simply use $W$. The minimizer $\pi^{\text{opt}}$ in \eqref{DefWassContinuous} is referred to as the \emph{optimal transportation plan}, and if $\mu\in\mathcal{P}_{2,{\rm{ac}}}(\mathbb{R}^{d})$, then  $\pi^{\text{opt}}$ is supported on the graph of the \emph{optimal transport map} $T^{\text{opt}}$ pushing $\mu_x$ to $\mu_y$. 

The connection between \eqref{AdditiveOptimizationMeasure} and \eqref{WassgradFlow} is that the minimizer of \eqref{AdditiveOptimizationMeasure} may be realized as the stationary solution of \eqref{WassgradFlow}. Conversely, if one is interested in the (possibly transient) solution of a PDE of the form \eqref{WassgradFlow}, then it might be possible to compute the same by performing discrete time-stepping realizing gradient descent for \eqref{AdditiveOptimizationMeasure}. 

In recent years, several algorithms have been proposed for solving measure-valued optimization problems, see e.g., \cite{benamou2016augmented,peyre2015entropic,carlier2017convergence,carrillo2021primal,mokrov2021large,alvarez2021optimizing}. In this work, we explore the possibility of leveraging the additive structure of the objective in \eqref{AdditiveOptimizationMeasure} for distributed nonparametric computation. 

\section{Main Idea}
We relabel the argument of the functional $F_{i}$ in \eqref{AdditiveOptimizationMeasure} as $\mu_{i}$ for all $i\in[n]$, and then impose the consensus constraint $\mu_{1} = \mu_{2} = \hdots = \mu_{n}$. Denoting $$\mathcal{P}_{2}^{n+1}(\mathbb{R}^{d}) := \underbrace{\mathcal{P}_{2}(\mathbb{R}^{d}) \times \hdots \times  \mathcal{P}_{2}(\mathbb{R}^{d})}_{n+1\;\text{times}},$$ 
we the rewrite \eqref{AdditiveOptimizationMeasure} as
{\small{\begin{subequations}
\begin{align}
&\underset{(\mu_1, \hdots, \mu_n, \zeta)\in \mathcal{P}_{2}^{n+1}(\mathbb{R}^{d})}{\arg\inf} \: F_{1}(\mu_1) + F_{2}(\mu_2) + \hdots + F_{n}(\mu_n) \label{AdditiveOptimizationnplusoneObj}\\
&\;\qquad\text{subject to} \qquad\; \mu_{i} = \zeta \quad\text{for all}\;i\in[n].\label{AdditiveOptimizationnplusoneConstr}
\end{align}
\label{AdditiveOptimizationnplusone}
\end{subequations}}}
Akin to the standard (Euclidean) augmented Lagrangian, we define the \emph{Wasserstein augmented Lagrangian} 
{\small{\begin{align}
&L_{\alpha}(\mu_1, \hdots, \mu_{n},\zeta,\nu_1,\hdots,\nu_{n}) := \nonumber\\
&\displaystyle\sum_{i=1}^{n}\bigg\{F_{i}(\mu_i) + \dfrac{\alpha}{2}W^{2}\left(\mu_i,\zeta\right) + \int_{\mathbb{R}^{d}}\nu_{i}(\bm{\theta})\left(\differential\mu_{i} -\differential\zeta\right) \bigg\}
\label{WassAugLagrangian}    
\end{align}}}
where $\nu_{i}(\bm{\theta})$, $i\in[n]$, are the Lagrange multipliers for the constraints in \eqref{AdditiveOptimizationnplusoneConstr}, and $\alpha>0$ is a regularization constant. 

Motivated by the Euclidean alternating direction method of multipliers (ADMM), we set up the recursions
{\small{\begin{subequations}
\begin{align}
\mu_{i}^{k+1} &= \underset{\mu_{i}\in\mathcal{P}_{2}(\mathbb{R}^{d})}{\arg\inf}\:L_{\alpha}\left(\mu_1, \hdots, \mu_{n},\zeta^{k},\nu_1^{k},\hdots,\nu_{n}^{k}\right) \label{ADMMrecursionsmuupdate}\\
\zeta^{k+1} &= \underset{\zeta\in\mathcal{P}_{2}(\mathbb{R}^{d})}{\arg\inf}\:L_{\alpha}\left(\mu_1^{k+1}, \hdots, \mu_{n}^{k+1},\zeta,\nu_1^{k},\hdots,\nu_{n}^{k}\right) \label{ADMMrecursionszetaupdate}\\
\nu_{i}^{k+1} &= \nu_i^{k} + \alpha\left(\mu_{i}^{k+1} - \zeta^{k+1}\right)\label{ADMMrecursionsnuupdate}
\end{align}
\label{ADMMrecursions}
\end{subequations}}}
where $i\in[n]$, and the recursion index $k\in\mathbb{N}_{0}$ (the set of whole numbers $\{0,1,2,\hdots\}$). We view \eqref{ADMMrecursionsmuupdate}-\eqref{ADMMrecursionszetaupdate} as primal updates, and \eqref{ADMMrecursionsnuupdate} as dual ascent. 

Let $\nu_{\text{sum}}^{k}(\bm{\theta}) := \displaystyle\sum_{i=1}^{n}\nu_{i}^{k}(\bm{\theta})$, $ k\in\mathbb{N}_{0}$. Substituting \eqref{WassAugLagrangian} in \eqref{ADMMrecursions}, dropping the terms independent of the decision variable in the respective $\arg\inf$, and re-scaling, the recursions \eqref{ADMMrecursions} simplify to
{\small{\begin{subequations}
\begin{align}
\mu_{i}^{k+1} &= \underset{\mu_{i}\in\mathcal{P}_{2}(\mathbb{R}^{d})}{\arg\inf}\: \dfrac{1}{2}W^{2}\left(\mu_{i},\zeta^{k}\right) \!+\! \dfrac{1}{\alpha}\bigg\{F_{i}(\mu_i) + \!\!\int_{\mathbb{R}^{d}}\!\! \nu_{i}^{k}(\bm{\theta})\differential\mu_{i}\bigg\}\nonumber\\
&= \prox^{W}_{\frac{1}{\alpha}\left(F_{i}(\cdot) + \!\int \nu_{i}^{k}\differential(\cdot)\right)}\left(\zeta^{k}\right), \label{ADMMrecursionsmuupdateSimplified}\\
\zeta^{k+1} &= \underset{\zeta\in\mathcal{P}_{2}(\mathbb{R}^{d})}{\arg\inf}\:\displaystyle\sum_{i=1}^{n}\bigg\{\dfrac{1}{2}W^{2}\left(\mu_{i}^{k+1},\zeta\right) - \dfrac{1}{\alpha}\!\int_{\mathbb{R}^{d}}\!\! \nu_{i}^{k}(\bm{\theta})\differential\zeta\bigg\}\nonumber\\
&= \underset{\zeta\in\mathcal{P}_{2}(\mathbb{R}^{d})}{\arg\inf}\bigg\{\left(\displaystyle\sum_{i=1}^{n}W^{2}\left(\mu_{i}^{k+1},\zeta\right)\right) - \dfrac{2}{\alpha}\!\int_{\mathbb{R}^{d}}\!\! \nu_{\text{sum}}^{k}(\bm{\theta})\differential\zeta\bigg\}, \label{ADMMrecursionszetaupdateSimplified}\\
\nu_{i}^{k+1} &= \nu_i^{k} + \alpha\left(\mu_{i}^{k+1} - \zeta^{k+1}\right),\label{ADMMrecursionsnuupdateSimplified}
\end{align}
\label{ADMMSimplified}
\end{subequations}}}
wherein we use the notation $\prox^{W}_{G(\cdot)}(\zeta)$ to denote the \emph{Wasserstein proximal operator} of the functional $G(\cdot)$, acting on $\zeta\in\mathcal{P}_{2}\left(\mathbb{R}^{d}\right)$, given by
\begin{align}
\prox^{W}_{G(\cdot)}(\zeta) := \underset{\mu\in\mathcal{P}_{2}\left(\mathbb{R}^{d}\right)}{\arg\inf}\:\dfrac{1}{2}W^{2}\left(\mu,\zeta\right) + G(\mu).
\label{defWassProx}
\end{align}
We can view \eqref{defWassProx} as a generalization of the finite dimensional Euclidean proximal operator
\begin{align}
\prox^{\|\cdot\|_{2}}_{g}(\bm{z}) := \underset{\bm{x}\in\mathbb{R}^{d}}{\arg\inf}\:\dfrac{1}{2}\|\bm{x}-\bm{z}\|_{2}^{2} + g(\bm{x}).   
\label{defEuclideanProx}    
\end{align}
We refer to \eqref{ADMMSimplified} as the \emph{Wasserstein consensus ADMM} -- the notion generalizes its finite dimensional Euclidean counterpart in the sense \eqref{ADMMrecursionsmuupdateSimplified}-\eqref{ADMMrecursionszetaupdateSimplified} are analogues of the so-called $x$ and $z$ updates, respectively; see e.g., \cite[Ch. 5.2.1]{parikh2014proximal}. However, important difference arises in \eqref{ADMMrecursionszetaupdateSimplified} compared to its Euclidean counterpart due to the sum of squares of Wasserstein distances. In the Euclidean case, the corresponding $z$ update can be analytically performed in terms of the arithmetic mean of the $x$ updates. While \eqref{ADMMrecursionszetaupdateSimplified} does involve a \emph{generalized mean} of the updates from \eqref{ADMMrecursionsmuupdateSimplified}, we now have \emph{Wasserstein barycentric proximal} of a linear functional. In other words, \eqref{ADMMrecursionszetaupdateSimplified} amounts to computing the Wasserstein barycenter of $n$ measures $\{\mu_{1}^{k+1},\hdots,\mu_{n}^{k+1}\}$ with a linear regularization involving $\nu_{\text{sum}}^{k}$.

The proximal updates \eqref{ADMMrecursionsmuupdateSimplified} are closely related to the Wasserstein gradient flows generated by the respective (scaled) free energy functionals $$\Phi_{i}(\mu_i) := F_{i}(\mu_i) + \int_{\mathbb{R}^{d}} \nu_{i}^{k}\differential\mu_i, \quad \mu_{i}\in\mathcal{P}_{2}(\mathbb{R}^{d}), \quad i\in[n].$$
Under mild assumptions on $\Phi_i$, as $1/\alpha \downarrow 0$, the sequence $\{\mu_{i}^{k}(\alpha)\}_{k\in\mathbb{N}_{0}}$ generated by the proximal updates \eqref{ADMMrecursionsmuupdateSimplified} converge to the measure-valued solution trajectory $\widetilde{\mu}_{i}(t,\cdot)$, $t\in[0,\infty)$, generated by the initial value problems (IVPs)
\begin{align}
\dfrac{\partial\widetilde{\mu}_{i}}{\partial t} = - \nabla^{W}\Phi_{i}\left(\widetilde{\mu}_{i}\right), \; \widetilde{\mu}_{i}(t=0,\cdot) = \widetilde{\mu}_{i}^{0}(\cdot), \; i\in[n].    
\label{WassGradFlow}    
\end{align}
Thus, in a rather generic setting, performing the proximal updates \eqref{ADMMrecursionsmuupdateSimplified} in parallel across the index $i\in[n]$, amounts to performing distributed time updates for the approximate transient solutions of the IVPs \eqref{WassGradFlow}. 

Important examples of $F_{i}$ include $\int V(\bm{\theta}) \differential\mu_{i}(\bm{\theta})$ (potential energy for some suitable advection potential $V$), $\beta^{-1}\int \log\mu_{i}(\bm{\theta}) \differential\mu_{i}(\bm{\theta})$ (internal energy with the ``inverse temperature'' parameter $\beta>0$), $\int_{\mathbb{R}^{2d}} U(\bm{\theta},\bm{\sigma}) \differential\mu_{i}(\bm{\theta}) \differential\mu_{i}(\bm{\sigma})$ (interaction energy for some symmetric positive definite interaction potential $U$).

To numerically realize the recursions \eqref{ADMMSimplified}, we consider a sequence of discrete probability distributions $\{\bm{\mu}_{1}^{k}, \hdots, \bm{\mu}_{n}^{k}, \bm{\zeta}^{k}\}$ indexed by $k\in\mathbb{N}_{0}$ where each distribution is a probability vector of length $N\times 1$, representative of the respective probability values at $N$ samples. Thus, for each fixed $k\in\mathbb{N}_{0}$, the tuple $$\left(\bm{\mu}_{1}^{k}, \hdots, \bm{\mu}_{n}^{k}, \bm{\zeta}^{k}\right)\in\underbrace{\Delta^{N-1} \times \hdots \times\Delta^{N-1}}_{n+1\;\text{times}} =: \left(\Delta^{N-1}\right)^{n+1}.$$
Likewise, for each fixed $k\in\mathbb{N}_{0}$, the Lagrange multipliers $\left(\bm{\nu}_{1}^{k},\hdots,\bm{\nu}_{n}^{k}\right)\in\mathbb{R}^{nN}$, and $\bm{\nu}^{k}_{\text{sum}}=\displaystyle\sum_{i=1}^{n}\bm{\nu}^{k}_{i}\in\mathbb{R}^{N}$.

Given probability vectors $\bm{\xi}, \bm{\eta}\in\Delta^{N-1}$, let
$\Pi_{N}\left(\bm{\xi},\bm{\eta}\right) := \{\bm{M}\in\mathbb{R}^{N\times N} \mid \bm{M} \geq \bm{0}\;\text{(elementwise)}, \allowbreak \bm{M}\bm{1} = \bm{\xi},\allowbreak \bm{M}^{\top}\bm{1} = \bm{\eta}\}$.
Also, let $\bm{C}\in\mathbb{R}^{N\times N}$ denote the squared Euclidean distance matrix for the sampled data $\{\bm{\theta}_{r}\}_{r\in[N]}$ in $\mathbb{R}^{d}$, i.e., the entries of the matrix $\bm{C}$ are $\bm{C}(i,j) := \|\bm{\theta}_{i}-\bm{\theta}_{j}\|_{2}^{2}$ for all $i,j\in[N]$.  

For each $i\in[n]$ and $k\in\mathbb{N}_{0}$, we write the discrete version of \eqref{ADMMSimplified} as
{\small{\begin{subequations}
\begin{align}
&\bm{\mu}_{i}^{k+1} = \prox^{W}_{\frac{1}{\alpha}\left(F_{i}(\bm{\mu}_{i}) + \langle\bm{\nu}_{i}^{k},\bm{\mu}_{i}\rangle\right)}\left(\bm{\zeta}^{k}\right)   \nonumber\\
&= \underset{\bm{\mu}_{i}\in\Delta^{N-1}}{\arg\inf}\bigg\{\underset{\bm{M}\in\Pi_{N}\left(\bm{\mu}_{i},\bm{\zeta}^{k}\right)}{\min}\frac{1}{2}\langle\bm{C},\bm{M}\rangle + \frac{1}{\alpha}\left(F_{i}(\bm{\mu}_{i}) + \langle\bm{\nu}_{i}^{k},\bm{\mu}_{i}\rangle\right)\bigg\}, \label{MuUpdateDiscrete}\\
&\bm{\zeta}^{k+1} = \underset{\bm{\zeta}\in\Delta^{N-1}}{\arg\inf} \bigg\{\! \left(\!\displaystyle\sum_{i=1}^{n}\!\underset{\bm{M}_{i}\in\Pi_{N}\left(\bm{\mu}_{i}^{k+1},\bm{\zeta}\right)}{\min}\frac{1}{2}\langle\bm{C},\bm{M}_{i}\rangle\!\right) - \frac{2}{\alpha}\langle\bm{\nu}^{k}_{\text{sum}},\bm{\zeta}\rangle\!\bigg\}, \label{ZetaUpdateDiscrete}\\
&\bm{\nu}_{i}^{k+1} = \bm{\nu}_{i}^{k} + \alpha\left(\bm{\mu}_{i}^{k+1} - \bm{\zeta}^{k+1}\!\right), \label{NuUpdateDiscrete}
\end{align}
\label{ADMMdiscrete}  
\end{subequations}}}
wherein \eqref{MuUpdateDiscrete}-\eqref{ZetaUpdateDiscrete} used the discrete version of the squared Wasserstein distance.

Replacing the squared Wasserstein distance in \eqref{ADMMSimplified} by the entropy a.k.a. Sinkhorn regularized squared Wasserstein distance, modify the recursions \eqref{ADMMdiscrete} as
\begin{subequations}
\begin{align}
&\bm{\mu}_{i}^{k+1} = \prox^{W_{\varepsilon}}_{\frac{1}{\alpha}\left(F_{i}(\bm{\mu}_{i}) + \langle\bm{\nu}_{i}^{k},\bm{\mu}_{i}\rangle\right)}\left(\bm{\zeta}^{k}\right)   \nonumber\\
&= \underset{\bm{\mu}_{i}\in\Delta^{N-1}}{\arg\inf}\bigg\{\underset{\bm{M}\in\Pi_{N}\left(\bm{\mu}_{i},\bm{\zeta}^{k}\right)}{\min}\bigg\langle\frac{1}{2}\bm{C} + \varepsilon\log\bm{M},\bm{M}\bigg\rangle \nonumber\\
&\qquad+ \frac{1}{\alpha}\left(F_{i}(\bm{\mu}_{i}) + \langle\bm{\nu}_{i}^{k},\bm{\mu}_{i}\rangle\right)\bigg\}, \label{MuUpdateDiscreteRegularized}\\
&\bm{\zeta}^{k+1} = \underset{\bm{\zeta}\in\Delta^{N-1}}{\arg\inf} \!\bigg\{\!\! \left(\!\displaystyle\sum_{i=1}^{n}\!\underset{\bm{M}_{i}\in\Pi_{N}\left(\bm{\mu}_{i}^{k+1},\bm{\zeta}\right)}{\min}\bigg\langle\frac{1}{2}\bm{C} + \varepsilon\log\bm{M}_{i},\bm{M}_{i}\!\bigg\rangle\!\right) \nonumber\\
&\qquad\qquad - \frac{2}{\alpha}\langle\bm{\nu}^{k}_{\text{sum}},\bm{\zeta}\rangle\!\bigg\}, \label{ZetaUpdateDiscreteRegularized}\\
&\bm{\nu}_{i}^{k+1} = \bm{\nu}_{i}^{k} + \alpha\left(\bm{\mu}_{i}^{k+1} - \bm{\zeta}^{k+1}\right), \label{NuUpdateDiscreteRegularized}
\end{align}
\label{ADMMdiscreteRegularized}  
\end{subequations}
where $\varepsilon > 0$ is a regularization parameter. In the remaining, we summarize novel results that enable us to numerically perform the recursions \eqref{ADMMdiscreteRegularized}.

\section{Results}
\subsection{The $\bm{\mu}$ Update}\label{subsec:MuUpdate}
The Sinkhorn regularized recursion \eqref{MuUpdateDiscreteRegularized} allows us to get semi-analytical handle on the nested minimization via strong duality. Specifically, consider the convex functions $F_i, G_i:\Delta^{N-1}\mapsto\mathbb{R}$ for all $i\in[n]$ where 
$G_{i}(\bm{\mu_i}) := F_{i}(\bm{\mu}_{i}) + \langle\bm{\nu}_{i}^{k},\bm{\mu}_{i}\rangle$, and denote the Legendre-Fenchel conjugate of $G_{i}$ as $G_{i}^{*}$. Following \cite[Lemma 3.5]{karlsson2017generalized}, \cite[Sec. III]{caluya2019gradient}, the Lagrange dual problem associated with \eqref{MuUpdateDiscreteRegularized} is
{\small{\begin{align}
&\left(\bm{\lambda}_{0i}^{\text{opt}},\bm{\lambda}_{1i}^{\text{opt}}\right) = \underset{\bm{\lambda}_{0i},\bm{\lambda}_{1i}\in\mathbb{R}^{N}}{\arg\max} \bigg\{\langle\bm{\lambda}_{0i},\bm{\zeta}_{k}\rangle - G_{i}^{*}\left(-\bm{\lambda}_{1i}\right) \nonumber\\
&- \alpha\varepsilon\left(\exp\left(\frac{\bm{\lambda}_{0i}^{\top}}{\alpha\varepsilon}\right)\exp\left(-\frac{\bm{C}}{2\varepsilon}\right)\exp\left(\frac{\bm{\lambda}_{1i}}{\alpha\varepsilon}\right)\right)\bigg\}, \; i\in[n].
\label{LagrangeDualProblem}    
\end{align}}}
Using \eqref{LagrangeDualProblem}, the proximal updates in \eqref{MuUpdateDiscreteRegularized} can then be recovered from the following proposition.
\begin{proposition}(\cite[Lemma 3.5]{karlsson2017generalized},\cite[Theorem 1]{caluya2019gradient})\label{prop:ProximalUpdateGeneral}
Given $\alpha,\varepsilon > 0$, the squared Euclidean distance matrix $\bm{C}\in\mathbb{R}^{N\times N}$, and the probability vector $\bm{\zeta}^{k}\in\Delta^{N-1}$, $k\in\mathbb{N}_{0}$. Let $\bm{0}$ denote the $N\times 1$ vector of zeros. For $i\in[n]$, the vectors $\bm{\lambda}_{0i}^{\rm{opt}},\bm{\lambda}_{1i}^{\rm{opt}}\in\mathbb{R}^{N}$ in \eqref{LagrangeDualProblem} solve the system
{\small{\begin{subequations}
\begin{align}
&\exp\left(\frac{\bm{\lambda}_{0i}^{\rm{opt}}}{\alpha\varepsilon}\right) \odot \left( \exp\left(-\frac{\bm{C}}{2\varepsilon}\right)\exp\left(\frac{\bm{\lambda}_{1i}^{\rm{opt}}}{\alpha\varepsilon}\right) \right) = \bm{\zeta}_{k},  \label{ZetakEquation}\\ 
&\bm{0} \in \partial_{\bm{\lambda}_{1i}^{\rm{opt}}}G_{i}^{*}\left(\!-\bm{\lambda}_{1i}^{\rm{opt}}\!\right) - \!\exp\!\left(\!\frac{\bm{\lambda}_{1i}^{\rm{opt}}}{\alpha\varepsilon}\!\right) \odot \left(\!\exp\!\left(\!-\frac{\bm{C}^{\top}}{2\varepsilon}\!\right)\!\exp\!\left(\!\frac{\bm{\lambda}_{0i}^{\rm{opt}}}{\alpha\varepsilon}\!\right)\!\right). \label{ZeroInSubdifferential}  
\end{align}
\label{lambda0lambda1equations}  
\end{subequations}}}
The proximal update $\bm{\mu}_{i}^{k+1}$ in \eqref{MuUpdateDiscreteRegularized} is given by
{\small{\begin{align}
\bm{\mu}_{i}^{k+1} = \exp\left(\frac{\bm{\lambda}_{1i}^{\rm{opt}}}{\alpha\varepsilon}\right) \odot \left(\exp\left(-\frac{\bm{C}^{\top}}{2\varepsilon}\right) \exp\left(\frac{\bm{\lambda}_{0i}^{\rm{opt}}}{\alpha\varepsilon}\right)\right).
\label{MuUpdateGeneral}    
\end{align}}}
\end{proposition}
We point out an important special case: if $F_{i}(\bm{\mu}_{i}) = \beta^{-1}\langle\log\bm{\mu}_{i},\bm{\mu}_{i}\rangle$ where $\beta>0$, then Proposition \ref{prop:ProximalUpdateGeneral} reduces exactly to \cite[Theorem 1]{caluya2019gradient} allowing further simplification of \eqref{ZeroInSubdifferential}. Then, \eqref{lambda0lambda1equations} can be solved via certain cone-preserving block coordinate iteration proposed in \cite[Sec. III.B,C]{caluya2019gradient} that is provably contractive. This makes the proximal update \eqref{MuUpdateGeneral} semi-analytical in the sense the pair $\left(\bm{\lambda}_{0i}^{\text{opt}},\bm{\lambda}_{1i}^{\text{opt}}\right)$ needs to be numerically computed by performing the block coordinate iteration while ``freezing'' the index $k\in\mathbb{N}_{0}$. With the converged pair $\left(\bm{\lambda}_{0i}^{\text{opt}},\bm{\lambda}_{1i}^{\text{opt}}\right)$, the evaluation \eqref{MuUpdateGeneral} is analytical for each $k\in\mathbb{N}_{0}$. 

In our context, another case of interest is when $F_{i}$ and hence $G_{i}$, is linear in $\bm{\mu}_{i}$. The following result shows that the proximal update $\bm{\mu}_{i}^{k+1}$ in this case can be computed analytically, obviating the zero order hold sub-iterations mentioned above.
\begin{theorem}\label{Thm:WassProxOfLinear}
Given $\bm{a}\in\mathbb{R}^{N}\setminus\{\bm{0}\}$, let $\bm{\Phi}(\bm{\mu}) := \langle\bm{a},\bm{\mu}\rangle$ for $\bm{\mu}\in\Delta^{N-1}$. Let $\bm{C}\in\mathbb{R}^{N\times N}$ be the squared Euclidean distance matrix, and for $\varepsilon>0$, let $\bm{\Gamma}:=\exp\left(-\bm{C}/2\varepsilon\right)$. For any $\bm{\zeta}\in\Delta^{N-1}$, $\alpha > 0$, we have
{\small{\begin{align}
\prox^{W_{\varepsilon}}_{\frac{1}{\alpha}\bm{\Phi}}\left(\bm{\zeta}\right) = \!\exp\!\left(\!-\dfrac{1}{\alpha\varepsilon}\bm{a}\!\right) \!\odot\! \left(\!\bm{\Gamma}^{\top}\!\left(\!\bm{\zeta}\!\oslash\!\left(\!\bm{\Gamma}\exp\!\left(\!-\dfrac{1}{\alpha\varepsilon}\bm{a}\!\right)\!\right)\!\right)\!\right).
\label{WassP}
\end{align}}}
\end{theorem}

\subsection{The $\bm{\zeta}$ Update}\label{subsec:ZetaUpdate}
The update \eqref{ZetaUpdateDiscreteRegularized} can be seen as a problem of computing the Sinkhorn regularized Wasserstein barycenter with an extra linear regularization. Let $W_{\varepsilon,\bm{\mu}_{i}}^{2}(\bm{\zeta}) := \underset{\bm{M}_{i}\in\Pi_{N}\left(\bm{\mu}_{i},\bm{\zeta}\right)}{\min}\bigg\langle\frac{1}{2}\bm{C} + \varepsilon\log\bm{M}_{i},\bm{M}_{i}\bigg\rangle$, $\varepsilon > 0$, for given $\bm{\mu}_{i}\in\Delta^{N-1}$ for all $i\in[n]$, and for a given squared Euclidean distance matrix $\bm{C}\in\mathbb{R}^{N\times N}$. Let the superscript $^{*}$ denote the Legendre-Fenchel conjugate. Following \cite[Sec. 4.1]{cuturi2016smoothed}, some calculations show that the dual problem corresponding to \eqref{ZetaUpdateDiscreteRegularized} becomes
{\small{\begin{align}
\left(\bm{u}_{1}^{\rm{opt}},\hdots,\bm{u}_{n}^{\rm{opt}}\right) = &\underset{\left(\bm{u}_{1},\hdots,\bm{u}_{n}\right)\in\mathbb{R}^{nN}}{\arg\min}\displaystyle\sum_{i=1}^{n}\left(W_{\varepsilon,\bm{\mu}_{i}^{k+1}}^{2}\right)^{*}\left(\bm{u}_{i}\right)\nonumber\\
&\quad{\rm{subject\;to}}\quad \displaystyle\sum_{i=1}^{n}\bm{u}_{i} = \frac{2}{\alpha}\bm{\nu}^{k}_{\text{sum}}.    
\label{DualOfSinkhornBaryWithLinReg}    
\end{align}}}
Consequently, the update \eqref{ZetaUpdateDiscreteRegularized} can be performed by first solving the problem \eqref{DualOfSinkhornBaryWithLinReg}, and then invoking the primal-dual relation $\bm{\zeta}^{\rm{opt}} = \nabla_{\bm{u}_{i}}\left(W_{\varepsilon,\bm{\mu}_{i}}^{2}\right)^{*}\left(\bm{u}_{i}^{\rm{opt}}\right)\in\Delta^{N-1}\:\forall i\in[n]$, at the minimizer of \eqref{DualOfSinkhornBaryWithLinReg}. It turns out that \eqref{DualOfSinkhornBaryWithLinReg} leads to an inner layer Euclidean ADMM whose structure allows efficient distributed computation.

The results summarized above lead to an overall algorithm realizing operator splitting for gradient flows in the manifold of probability measures, which solve \eqref{AdditiveOptimizationMeasure} via distributed computation. Numerical experiments (not reported herein due to page constraints) on several test problems of the form \eqref{AdditiveOptimizationMeasure} reveal that the proposed framework has good computational performance.

\vspace*{-0.05in}

\bibliography{ifacconf}             
                                                   







\end{document}